\documentclass[12pt,reqno]{amsart}

\usepackage{hyperref, amssymb, amsmath, enumerate, verbatim, color, bbm}
\usepackage[dvipsnames]{xcolor}

\advance \topmargin by -\headheight
\advance \topmargin by -\headsep

\evensidemargin \oddsidemargin
\hypersetup{pdfstartview={FitH}}

\newtheorem{thm}{Theorem}[section]
\newtheorem{prop}[thm]{Proposition}

\newtheorem*{prop*}{Proposition}

\newtheorem{lem}[thm]{Lemma}
\newtheorem{lemma}[thm]{Lemma}

\newtheorem{conj}[thm]{Conjecture}

\theoremstyle{definition}
\newtheorem{definition}[thm]{Definition}

\theoremstyle{remark}
\newtheorem{remark}[thm]{Remark}
\numberwithin{equation}{section}
\newtheorem{claim}{Claim}

\newcommand{\NN}{\mathbb{N}}
\newcommand{\RR}{\mathbb{R}}

\newcommand{\FF}{\mathcal{F}}

\newcommand{\TT}{\mathbb{T}}
\newcommand{\FQ}{\mathbb{F}_q}
\newcommand{\ZZ}{\mathbb{Z}}

\newcommand{\PP}{\mathbb{P}}

\usepackage{amsmath, amsthm, amssymb, stmaryrd}
\usepackage{hyperref}
\usepackage{enumerate}
\usepackage{url}

\let\OLDthebibliography\thebibliography
\renewcommand\thebibliography[1]{
  \OLDthebibliography{#1}
  \setlength{\parskip}{0pt}
  \setlength{\itemsep}{0pt plus 0.3ex}
}

\usepackage{mathtools}

\usepackage{multirow, array}
\usepackage{placeins}

\usepackage{caption} 
\newcommand*\wrapletters[1]{\wr@pletters#1\@nil}
\def\wr@pletters#1#2\@nil{#1\allowbreak\if&#2&\else\wr@pletters#2\@nil\fi}

\def\alp{{\alpha}}

\def\del{{\delta}}

\def\kap{{\kappa}}
\def\eps{\varepsilon}

\def\le{\leqslant} \def\ge{\geqslant}  
\def \leq {\leqslant} \def \geq {\geqslant}

\def \lam {{\lambda}}

\def \bF {\mathbb F}

\def \bK {\mathbb K}

\def \bP {\mathbb P}

\def \bR {\mathbb R}
\def \bZ {\mathbb Z}
\def \bT {\mathbb T}

\def \cC {\mathcal C}

\def \cF {\mathcal F}

\def \cP {\mathcal P}

\def \cS {\mathcal S}

\def \dim {\mathrm{dim}}
\def \ord {\mathrm{ord}}

\def \deg {\mathrm{deg}}

\newenvironment{proofof}[1]{\indent{\scshape Proof of #1}:~~}{\qed}
\newenvironment{proofcap}{\indent{\scshape Proof}:~~}{\qed}

\begin{document}

\title{Lonely runners in function fields}

\author{Sam Chow}
\address{Sam Chow, Mathematical Institute, University of Oxford, Woodstock Road, Oxford OX2 6GG, United Kingdom}
\email{Sam.Chow@maths.ox.ac.uk}

\author{Luka Rimani\'{c}}
\address{Luka Rimani\'{c}, School of Mathematics, University of Bristol, University Walk, Bristol BS8 1TW, United Kingdom}
\email{luka.rimanic@bristol.ac.uk}


\begin{abstract}
The lonely runner conjecture, now over fifty years old, concerns the following problem. On a unit length circular track, consider $m$ runners starting at the same time and place, each runner having a different constant speed. The conjecture asserts that each runner is \emph{lonely} at some point in time, meaning distance at least $1/m$ from the others. We formulate a function field analogue, and give a positive answer in some cases in the new setting. 
\end{abstract}

\maketitle


\section{Introduction}
{Introduced by Wills \cite{Wi67} in 1967, and then independently by Cusick \cite{Cu73}, the \textit{lonely runner conjecture} (LRC)} concerns the following problem. On a unit length circular track one considers $m$ runners who all start at the same place and at the same time, each runner having a constant speed, with speeds being pairwise distinct. We say that a runner is \textit{lonely} if all the other runners are of distance at least $1/m$ from her. A well-known reformulation of the problem \cite{BHK01}, provided below, states the problem in terms of integer speeds, with one of the runners having zero speed. Throughout, let $\|x\|$ denote the distance from a real number $x$ to the nearest integer. 

\begin{conj}[LRC]\label{thm:lrcconj}
Let $D$ be a set of $k$ positive integers. Then there exists $t \in \RR$ such that 
\begin{equation*}
\| t \cdot d \| \geq \frac{1}{k+1}, \hspace{3mm} \forall d \in D.
\end{equation*}
\end{conj}

At present, the problem is open for $k\geq 7$, that is, for eight or more runners. For $k=1$ the conjecture is trivial, and the $k=2$ case is resolved during the first lap of the slower non-stationary runner. The case $k=3$ was solved by Betke and Wills \cite{BW72}, and Cusick \cite{Cu73, Cu74, Cu82} as a problem in diophantine approximation. Cusick and Pomerance \cite{CP84} established the conjecture for $k=4$ by extending Cusick's previous work with additional estimates on certain exponential sums, though their proof needed a computer check for certain cases. Later Biennia et al. \cite{BGGST98} presented a much simpler proof of the $k=3,4$ cases, and additionally connected the lonely runner conjecture with a conjecture on flows in matroids. The case $k=5$ was first solved by Bohman, Holzman and Kleitman \cite{BHK01}, and then Renault \cite{Re04} provided a shorter proof. Finally, Barajas and Serra \cite{BS08} settled the case $k=6$. 

There has been recent progress for large values of $k$. For instance, Dubickas \cite{Du11} employed the Lov\'asz local lemma \cite{AS08} to establish the conjecture for certain lacunary sequences of speeds. For any set $D$ of $k$ non-zero integers, define
\begin{equation*}
\delta_{k} (D) := \sup_{t \in \RR /\ZZ} \min_{d\in D} \|dt\|,
\end{equation*}
and let $\delta_k$ be the infimum over all such sets $D$. We note that Dirichlet's approximation theorem in the form \cite[Lemma 2.1]{Va97} yields $\delta_k \leq \frac{1}{k+1}$, and the lonely runner conjecture states that $\delta_k \geq \frac{1}{k+1}$. The union bound 
\[
\bP_{t \in \bR / \bZ} \Bigl(\| dt \| \ge \frac1{2k + \eps}, \hspace{3mm} \forall d \in D \Bigr) \ge 1 - k \frac2{2k+\eps} > 0
\]
implies that $\del_k \ge \frac{1}{2k}$. Until recently, previous work \cite{Ch94, CC99, PS16} had only improved the denominator by an additive constant. Tao \cite{Ta17} showed that there exists an absolute constant $c > 0$ such that
\begin{equation}\label{tao_large}
\delta_k \geq \frac{1}{2k} + \frac{c\log k}{k^2 (\log \log k)^2}
\end{equation} 
for all sufficiently large $k$. In addition, he showed that it is enough to verify the conjecture for speeds at most $k^{O(k^2)}$, and also that the conjecture holds when the speeds are all small.
\begin{thm}[\cite{Ta17}]\label{tao_small}
Let $k \geq 1$ and suppose that $\max_{d\in D} |d| \leq 1.2k$. Then $\delta_k (D) \geq \frac{1}{k+1}$.
\end{thm}

Czerwi\' nski and Grytczuk \cite{CG08} proved that the conjecture holds provided one is allowed to make a runner invisible to the other runners.
\begin{thm}[\cite{CG08}]\label{invisible}
Let $k$ and $s$ be integers such that $0 \leq s <k$. Then for every set $D$ of $k$ positive integers there exists a subset $S\subseteq D$ of size $k-s$ such that 
\[
\delta_{k-s} (S) \geq \frac{s+1}{2k}.
\]
\end{thm}
\noindent Recently Perarnau and Serra \cite{PS16} extended this result by showing that there exists a time when either a runner is lonely or four runners are ``almost lonely''. Finally, Czerwi\'nski  \cite{Cz12} confirmed the conjecture for a random set of positive integer speeds (see also \cite{Al13}), establishing the inequality 
\[
\delta_{k} (D) \ge \frac12 - \eps
\]
with high probability, for any $\eps > 0$.
\\

In this article, we formulate an analogue of the lonely runner problem in function fields, and prove results of a similar flavour to (\ref{tao_large}), Theorem \ref{tao_small} and Theorem \ref{invisible}. We shall see, however, that our methods are very different.

\subsection{Function fields}

The analogy between number fields and function fields \cite{vdGMS05} has a long and distinguished history, dating back at least to a famous 1882 paper by Weber and Dedekind \cite{WD1882}. A significant milestone was reached by Weil \cite{We41},  who established the Riemann hypothesis for algebraic curves defined over finite fields, otherwise known as the ``Riemann hypothesis for function fields''. Since then, this connection has been deeply investigated in many contexts, and function field models have provided a valuable testing ground for building intuition, making predictions and developing proof techniques \cite{EVW16, En16, GG17, KS1, KS2, LW10}.\\

Let $q$ be a prime power, and let $\FQ$ denote the field of $q$ elements. We have the following analogy---see for instance \cite[\S2]{LW10}. 
\begin{center}
\begin{tabular}{ |c|c| } 
 \hline
 $\ZZ$ & $\FQ [T]$  \\ 
 $\RR$ & $\bK_\infty := \FQ ((T^{-1}))$ \\ 
 $x \in \RR / \ZZ$ & $x = \sum_{i = 1}^\infty x_{-i} T^{-i}$\\ 
 \hline
\end{tabular}
\end{center}
In order to measure loneliness, we require a notion of distance. For
\[
\alp = \sum_{i = - \infty}^n x_i T^i \in \bF_q((T^{-1})),
\]
let $\ord(\alp)$ be the greatest integer $i$ for which $x_i \ne 0$, and write $\langle \alp \rangle = q^{\ord (\alp)}$. We adopt the convention that $\ord (0) = -\infty$ and $\langle 0 \rangle = 0$. Our analogue of the unit track shall be the compact additive subgroup
\[
\bT =  \{\alp \in \bK_\infty: \langle \alp \rangle < 1 \}
\]
of $\bK_\infty$, and the ``norms" in $\TT$ take values in $\{ 0 \} \cup \{ q^{-1}, q^{-2}, \ldots \}.$
Observe that any element $\alp \in \bK_\infty$ can be uniquely decomposed as
\[
\alp = [\alp] + \| \alp \|,
\]
where $[\alp] \in \bF_q[t]$ and $\|\alp\| \in \bT$; this corresponds to decomposing an integer into its integer and fractional parts. For convenience, write
\[
|\alp| := \langle \| \alp \| \rangle;
\]
this corresponds to the distance from a real number to the nearest integer.

In analogy with the notation from the real setting, we define the \emph{loneliness} of a set $\cF \subseteq \bF_q[T] \setminus \{ 0 \}$ by
\[
\delta (\cF) := \sup_{\alp \in \bT} \min_{f \in \cF} |\alp f|.
\]
We formulate the lonely runner conjecture in the function field setting as follows.

\begin{conj}[LRC in function fields]\label{thm:lrcffconj}
Let $\cF \subseteq \bF_q[T] \setminus \{0\}$ be such that 
\[
1\leq |\cF| < \frac{q^{k+1}-1}{q-1}.
\]
Then
\[
\delta (\cF) \ge q^{-k}.
\]
\end{conj}

We shall assume throughout that all polynomials in $\cF$ are monic. We lose nothing in doing so, for if $c \in \bF_q \setminus \{ 0 \}$ and $f \in \bF_q[T]$ then $|\alp f| = |\alp(cf)|$.

\subsection{Principal findings}

The upper bound on $|\cF|$ is necessary, as illustrated by the following example. Let 
\begin{equation}\label{counterexample}
\FF_k := \bigcup_{j=0}^{k} \left\{ T^j + i_{j-1} T^{j-1} + \cdots + i_1 T + i_0 \colon i_0, \ldots, i_{j-1} \in \FQ \right\}.
\end{equation}
Note that $|\FF_k| = \sum_{j=0}^{k} q^{j} = \frac{q^{k+1}-1}{q-1}$, and we claim that $\delta(\FF_k)\leq q^{-(k+1)}$. Indeed, for $\alpha \in \TT$, the system of $k$ linear equations in $k+1$ variables
\begin{equation*}
\langle \alpha f \rangle [T^{-1}] = \langle \alpha f \rangle [T^{-2}] = \cdots = \langle \alpha f \rangle [T^{-k}]  = 0
\end{equation*}
always has solutions in $\bF_{q}^{k+1}$, where $\langle g \rangle [T^{m}]$ denotes the coefficient of $T^m$ in $g$. Therefore, for each $\alpha \in \TT$ we can find $f\in\FF_k$ such that $|\alpha f| \le q^{-(k+1)}$, verifying our claim.

The lower bounds on the size threshold turn to be much more demanding, and we provide only partial answers in this article. Our first result is that $q^k$ is the correct order of magnitude. The theorem below confirms Conjecture \ref{thm:lrcffconj} for $k=1$.

\begin{thm}\label{maximalloneliness}
Let $\cF \subseteq \bF_q[T] \setminus \{0\}$ be of size at most $q^k$. Then $\delta(\cF) \geq q^{-k}$.
\end{thm}

By essentially the same proof, one finds that maximal loneliness is always attained when the coefficient field is infinite. Precisely, if $K$ is an infinite field and $\cF$ is a finite subset of $K[T] \setminus \{ 0 \}$ then there exists 
\[
\alp = \sum_{i = 1}^\infty x_{-i} T^{-i} \qquad (x_{-i} \in K)
\]
such that $(\alp f) [T^{-1}] \ne 0$, for all $f \in \cF$.\\

Our second result is that for $k>1$, Conjecture \ref{thm:lrcffconj} holds provided that $D := \displaystyle \max_{f\in\FF}\deg (f) \le q - O(1)$, which can be thought of as an analogue of Theorem \ref{tao_small}. In order to state the result we denote by $N(m,q)$ the number of irreducible monic polynomials of degree $m$ over $\FQ$. Gauss established the well-known formula
\[
N(m,q) = \frac{1}{m} \sum_{d\mid m} \mu (d) q^{m/d},
\]
where $\mu$ is the \textit{M\"obius function}. In particular
\[
N(m,q) \ge \frac{q^m - q^{m-1} - \cdots - q}m.
\]

\begin{thm}[Small degrees]\label{smallD}
Let $k > 1$, let $q$ and $D$ be such that 
\[
\frac{N(k+1, q)}{q^{k} + \cdots + q} > \Bigl\lfloor \frac{D}{k+1} \Bigr\rfloor .
\]
Then for any set of non-zero polynomials $\FF$ of size at most $\frac{q^{k+1}-1}{q-1} - 1$ whose degrees are at most $D$, one has $\delta (\cF) \geq q^{-k}$. 
\end{thm}

Our final result is a non-trivial lower bound on the number of polynomials, irrespective of $D$. We are able to get close to a halfway between Theorem \ref{maximalloneliness} and Conjecture \ref{thm:lrcffconj} in the case $k=2$. 

\begin{thm}[Almost halfway there]\label{close_half}
There exists a universal constant $C$ such that the following holds. Let $\cF$ be a set of non-zero polynomials such that
\[
|\FF| \leq q^2 + 0.4877q - C.
\]
Then $\cF$ is of loneliness at least $q^{-2}$.
\end{thm}
\noindent Carefully following the proof of Theorem \ref{close_half}, one could obtain $C$ to be $6$.

\subsection{Methods} 

The first step is to recast the question at hand as a covering problem. In the process we associate to each runner a \emph{partial circulant matrix}; these are well-studied \cite{Dav1994} and, for instance, have important applications to compressed sensing \cite{RRT12, YMYZ10}. Once viewed as a covering problem, Theorem \ref{maximalloneliness} follows straightforwardly from the union bound. We prove Theorem \ref{smallD} by exploiting the structure that irreducible polynomials bring to the covering problem. We argue that if the set of speeds lacks multiples of an irreducible polynomial of degree $k+1$, then one requires many runners to cover the missing structure. 

Our proof of Theorem \ref{close_half} is considerably more involved. It uses a very particular structure known as a \emph{sunflower}---see for instance \cite{Ju11}. Formally, a sunflower is a family of subsets for which there exists a \textit{core} $K$ such that for every pair of distinct subsets in the family, their intersection is precisely $K$. Large sunflowers whose core has generic dimension are very efficient at covering many vectors, and these play an indispensable role in our proof. Several authors have investigated upper bounds on the size of a sunflower-free family of sets, and this topic has received some attention recently \cite{NS2016, He2017}. On the other hand, sunflowers have a fascinating structure, and there appear to be few instances in which this structure has been brought to bear on a separate problem. In the proof of Theorem \ref{close_half} we are able to exploit the structure inherent in vector subspaces forming a sunflower. Sunflowers in this linear setting have been useful in coding theory, see \cite{ER2015} and the references within. We consider our use of sunflowers to be an interesting feature in its own right. We ultimately consider two cases, according to whether or not the set of speeds contains a large sunflower.

\subsection{Organisation} In \S \ref{geometric_int} we present our covering interpretation, and use it to establish Theorem \ref{maximalloneliness}. Theorems \ref{smallD} and \ref{close_half} are proved in Sections \ref{small_D_proof} and \ref{halfway_section} respectively.

\subsection{Acknowledgements}

The first named author was supported by EPSRC Programme Grant EP/J018260/1 at the University of York, and by the National Science Foundation under Grant No. DMS-1440140 while in residence at the Mathematical Sciences Research Institute in Berkeley, California, during the Spring 2017 semester. The second named author was supported by the School of Mathematics at the University of Bristol. This work began when the authors were graduate students together at the University of Bristol, and forms part of the second named author's dissertation. We thank Julia Wolf and Trevor Wooley for their enthusiastic guidance and unwavering support. In particular, we are grateful to the latter for suggesting the project. Thanks to Tom Bloom and the rest of the HARICOT members for useful discussions. Finally, we are grateful to the anonymous referees for their feedback.

\section{A covering problem}\label{geometric_int}
In this section we connect Conjecture \ref{thm:lrcffconj} with a covering problem in which we cover a vector space by certain subspaces. We will be studying the lonely runner conjecture in this formulation for the remainder of the paper. We finish the section by proving Theorem \ref{maximalloneliness}.

We will use the following setup for the rest of the paper, unless stated otherwise. For $k, D \in \NN$ let $\cF$ be a set of non-zero polynomials in $\FQ [T]$ of maximal degree at most $D$, and suppose we want to check whether $\delta(\cF) \geq q^{-k}$. We can write any $f \in \cF$ as 
\[
f = a_{0}^{(f)} + a_{1}^{(f)}T + \cdots + a_{D}^{(f)}T^{D},
\]
for some non-zero vector $(a_0^{(f)},\ldots,a_D^{(f)}) \in \bF_q^{D+1}$. Note that there is a one-to-one correspondence between vectors in $\FQ^{D+k}$ and polynomials of degree at most $D+k-1$. To each $f \in \cF$ we associate a partial circulant matrix
\[
A_f = \begin{bmatrix}
    a_{0}^{(f)}       & a_{1}^{(f)}  & \dots & a_{k-1}^{(f)} & \dots & a_{D}^{(f)} &  &  & \\
      & a_{0}^{(f)}  & \dots & a_{k-2}^{(f)} & \dots & a_{D-1}^{(f)} & a_{D}^{(f)}  &  &  \\
     &  & \ddots & \ldots & \ldots & \ldots & \ldots & \ddots &  \\
      &   &  & a_{0}^{(f)} & \dots & a_{D-k+1}^{(f)} & a_{D-k+2}^{(f)}& \dots & a_{D}^{(f)}
\end{bmatrix} = \begin{bmatrix}
    f \\
      Tf  \\
     \vdots   \\
      T^{k-1}f
\end{bmatrix},
\] 
whose size is $k\times (D+k)$. We may assume that a time $\alpha$ takes the form
\[
\alpha = x_{-1}T^{-1} + x_{-2}T^{-2} + \cdots + x_{-(D+k)} T^{-(D+k)}.
\]
since further terms do not affect whether or not $|\alpha f| \geq q^{-k}$. We also associate to $\alp$ the column vector
\[
x= [x_{-1} \ x_{-2} \ \ldots \ x_{-(D+k)} ]^T.
\]
Now $|\alpha f | \geq q^{-k}$ for all $f \in \cF$ is equivalent to
\[
A_f x \neq 0, \hspace{3mm} \forall f \in \cF.
\] 
It is easy to see that the vector subspace 
\[
\ker (A_f)  = \langle f, Tf, \ldots, T^{k-1}f \rangle^{\perp} \leq \FQ^{D+k}
\]
has codimension $k$. Defining $\ker (f) = \ker (A_f)$, we obtain the following.

\begin{lem}[Covering version]\label{lrc_covering}
A set $\cF$, comprising non-zero polynomials of degree at most $D$, is of loneliness at least $q^{-k}$ if and only if
\[
\bigcup_{f \in \cF}\ker (f) \neq \FQ^{D+k}.
\]
\end{lem}

This simple reformulation of the problem already provides us with the means to prove Theorem \ref{maximalloneliness}. One only needs to observe that each subspace contains the origin, and apply the union bound.

\begin{proofof}{Theorem \ref{maximalloneliness}}
As mentioned, for each $f \in \cF$, the subspace $\ker (f)$ has codimension $k$ in $\FQ^{D+k}$, meaning that each $\ker (f)$ covers at most $q^{D}-1$ non-zero vectors of $\FQ^{D+k}$. Therefore, by the union bound, assuming $|\cF| \leq q^{k}$, the number of non-zero vectors covered by $\bigcup_{f \in \cF} \ker (f)$ is at most
\[
\left| \bigcup_{f \in \cF} \ker (f) \setminus \{0\} \right| \leq \sum_{f \in \cF} \left| \ker (f) \setminus \{0 \} \right| \leq |\cF| (q^{D}-1) \leq q^{D+k} - q^k < q^{D+k} - 1,
\]
which gives $\bigcup_{f \in \cF} \ker (f) \neq \FQ^{D+k}$, from which Lemma \ref{lrc_covering} yields $\delta (\cF) \geq q^{-k}$.
\end{proofof}

\section{Low degree polynomials}\label{small_D_proof}

In this section we prove Theorem \ref{smallD}. To do so, we exploit the fact that irreducible polynomials have a particular kernel structure. In particular, we use this to prove that every irreducible polynomial of degree $k+1$ has to be observed in our set of speeds, at least as a factor. We formalise this in the following way.

\begin{lem}[Irreducible polynomials are factors]\label{irreducibles_k}
Let $\FF$ be a set of non-zero polynomials such that $|\FF| \leq \frac{q^{k+1}-1}{q-1}-1$ and $\displaystyle \cup_{f\in\FF} \ker (f) = \FQ^{D+k}$. Then for each irreducible monic polynomial $G_{k+1}$ of degree $k+1$, there exists $f\in \FF$ such that $G_{k+1} \mid f$.
\end{lem}

\begin{proofof}{Theorem \ref{smallD}, assuming Lemma \ref{irreducibles_k}}
For the sake of contradiction suppose that $\delta (\cF) < q^{-k}$, or equivalently, by Lemma \ref{lrc_covering}, that $\cup_{f\in \FF} \ker (f) = \FQ^{D+k}$. In view of the assumption $|\FF| \leq \frac{q^{k+1}-1}{q-1} -1$, Lemma \ref{irreducibles_k} tells us that for each irreducible polynomial $G_{k+1}$ of degree $k+1$ we can find $f\in\FF$ such that $G_{k+1} \mid f$. Since the degree of each $f\in\FF$ is bounded by $D$, a single $f$ can contain at most $\lfloor D/(k+1) \rfloor$ such $G_{k+1}$ as factors. Consequently, we need at least 
\[
\frac{N(k+1, q)}{\Bigl\lfloor \frac{D}{k+1} \Bigr\rfloor}
\]
polynomials in $\FF$. This yields 
\[
\frac{N(k+1, q)}{\Bigl\lfloor \frac{D}{k+1} \Bigr\rfloor} \leq |\FF| \leq \frac{q^{k+1}-1}{q-1} - 1 = q^k + \cdots + q,
\]
contradicting our hypotheses. This completes the proof of Theorem \ref{smallD}, given Lemma \ref{irreducibles_k}.
\end{proofof}

It remains to prove Lemma \ref{irreducibles_k}. The pivotal idea is that if an irreducible polynomial is not a divisor of some $f\in\FF$, we can lower bound the number of lines one needs to cover the missing $k+1$ dimensional space.

\begin{lem}[Covering with lines]\label{lines}
Let $V$ be a vector space over $\FQ$ such that $\dim (V) =k+1$, and let $V_1, \ldots, V_{R} \leq V$ be subspaces such that $\dim (V_i) =1$, $1\leq i \leq R$. If $R \leq \frac{q^{k+1}-1}{q-1} - 1$ then 
\[
V \neq \bigcup_{i=1}^{R} V_i.
\]
\end{lem}
\begin{proofcap}
We need to cover $q^{k+1} - 1$ non-zero vectors with $R$ lines through the origin. Note that each such line has at most $q-1$ non-zero points, so the union bound gives
\[
\Biggl| \bigcup_{i=1}^{R} V_i \setminus \{ 0\} \Biggr| \leq \sum_{i=1}^{R} |V_i \setminus \{ 0\}| \leq \left(  \frac{q^{k+1}-1}{q-1} - 1 \right) (q - 1) = q^{k+1} - q < q^{k+1} -1,
\]
concluding the proof of Lemma \ref{lines}.
\end{proofcap}

Let $\cP_{=m}$ denote the set of all monic polynomials of degree $m$, and write $\cP_{\leq m} = \cup_{i=0}^{m} \cP_{i}$. Define 
\[
\cP_{\leq m} \FF = \{ P f \colon P\in \cP_{\leq m}, f\in \FF \}.
\] 
Recall that in our covering problem we seek to prove that $\cup_{f\in\FF}\ker(f) \neq \FQ^{D+k}$. The following lemma provides a sufficient condition for this to hold.

\begin{lem}\label{intersection_nonempty}
Let $V \leq \FQ^{D+k}$ be such that $\dim (V^\perp) = k+1$, and let $\FF$ be such that $|\FF| \leq \frac{q^{k+1}-1}{q-1} - 1$. If $\cP_{\leq k-1} \FF \cap V = \emptyset$ then 
\[
\bigcup_{f\in \FF} \ker (f) \neq \FQ^{D+k}.
\]
\end{lem}
\begin{proofcap}
Let $V = \langle v_1, \ldots, v_{D-1} \rangle$. The assumption $\cP_{\leq k-1} \FF \cap V = \emptyset$ implies that for each $f\in \FF$ the set $B_f = \{ v_1, \ldots, v_{D-1}, f, Tf, \ldots, T^{k-1}f \}$ is a linearly independent set of vectors in $\FQ^{D+k}$, and so
\[
\dim (\text{span} (B_{f})^{\perp}) = 1, \hspace{3mm} \forall f \in \cF.
\]
Lemma \ref{lines} finishes the proof, since $|\FF| \leq \frac{q^{k+1}-1}{q-1} - 1$. Indeed, our kernels fail to cover $V^\perp$ so they cannot cover the entire space.
\end{proofcap}
\\

We conclude this section by proving Lemma \ref{irreducibles_k}.

\begin{proofof}{Lemma \ref{irreducibles_k}}
Define $V_{G_{k+1}} = \langle G_{k+1}, TG_{k+1}, \ldots, T^{D-2}G_{k+1} \rangle$. The assumption that $\cup_{f\in \FF} \ker (f) = \FQ^{D+k}$, combined with Lemma \ref{intersection_nonempty}, gives that there exist $P\in \cP_{\leq k-1}$ and $f\in\FF$ such that $Pf \in V_{G_{k+1}}$. Thus, there exists a polynomial $p$ such that $Pf = G_{k+1}p$. Now the fact that $G_{k+1}$ is an irreducible monic polynomial of degree $k+1$, coupled with the fact that $P$ is of degree at most $k-1$, yields $G_{k+1}\mid f$.
\end{proofof}


\section{Getting close to halfway when $k=2$}\label{halfway_section}
In this section we prove Theorem \ref{close_half}, making use of the combinatorial notion of a \textit{sunflower}. Although we do not use any of the existing results on sunflowers, we are confident that this language helps the reader to understand the main ideas, and that the questions arising are of independent interest. 

\begin{definition}
A collection of sets $\cS$ forms a \textit{sunflower} if there exists a set $K$ such that for each $S_1,S_2 \in \cS$ with $S_1 \ne S_2$, one has
\[
S_1 \cap S_2 = K.
\]
We say that $K$ is the \textit{core} and we call an element $S\in \cS$ a \textit{petal}. If $\cS$ is a collection of vector subspaces then we say that $\cS$ is a \textit{sunflower of codimension }$d$ if $K$ is of codimension $d$.
\end{definition}

When considering loneliness $\delta (\cF) \geq q^{-2}$, the ambient space is $\FQ^{D+2}$, and we will be interested in collections that are formed by the kernels of polynomials in $\FF$. A petal thus takes the form $\ker(f)$, for some polynomial $f \in \FF$. The fact that we work with subspaces provides us with an additional tool: we have a notion of dimension, which forces the cardinality of a subspace to be a power of $q$. Observe that a ``generic'' intersection of two kernels has codimension $4$, being the null space of a rank four matrix. In light of this, one might attempt to efficiently cover the space by starting with a large sunflower of codimension 4. This motivates us to consider a largest such sunflower, that is, one with the most petals.

We consider two cases: the first is when there exists a large codimension 4 sunflower in the set of speeds, and the second is when all codimension 4 sunflowers are small. We remark that in the former scenario it is possible to sharpen the bound provided here to get the full statement of Conjecture \ref{thm:lrcffconj}, see Proposition \ref{largeimpliesdone} and Remark \ref{remark_large} below. In this section, we prove a slightly weaker statement in order to make the paper easier to follow, as the real improvement needs to be made in the second case. We discuss the latter further in the appendix.

We divide the section into several parts, emphasising the connection between sunflowers and the covering statement in Lemma \ref{lrc_covering}. We think of $\ker (f)$, $f\in \cF$, as being introduced in some order. In order to improve on the union bound, we shall consider the contribution of each polynomial with respect to previous ones, meaning the number of points that $\ker (f)$ covers that were not already covered by the previously introduced polynomials.\\

The objective of our first lemma is to explain how the existence of a sunflower in the set of speeds reduces the contribution of the remaining polynomials. Moreover, it highlights the role played by the size of a largest sunflower.

\begin{lem}[After the sunflower]\label{afterparty}
Let $\cS = \{ f_1,\ldots,f_n\}$ be a sunflower of codimension $4$ of maximal size in $\cF$, and let $f\in\cF \setminus \cS$. Then the contribution of $f$ is bounded from above by
\begin{equation}\label{after_contribution}
\left| \ker(f) \setminus \bigcup_{i=1}^{n} \ker (f_i) \right| \leq \max \{ q^D - q^{D-1}, \: \: q^{D} - n q^{D-2} + (n-1) q^{D-3} \}.
\end{equation}
\end{lem}
\begin{proofcap}
To ease notation define $K_i = \ker (f_i)$, for $1\leq i \leq n$, and write $K_f = \ker (f)$. Let $K$ be the core of $\cS$, meaning that for all $1\leq i < j \leq n$ we have $K = K_i \cap K_j$. The contribution from $f$ is bounded above by
\[
\left| K_f \setminus \bigcup_{i=1}^{n} K_i \right| = q^{D} - \left| V_1 \cup \cdots \cup V_n \right|,
\]
where $V_i = K_i \cap K_{f}$ ($1\leq i \leq n$). Note that 
\begin{align*}
\left| V_1 \cup \cdots \cup V_n \right|  &=  |V_1| +	 \sum_{j=2}^{n} \left| V_j \setminus \left( V_1 \cup \cdots \cup V_{j-1} \right) \right| \\
&= |V_1| + \sum_{j=2}^{n} \left( |V_j| - \left| (V_j \cap V_1) \cup \cdots \cup (V_j \cap V_{j-1}) \right| \right).
\end{align*}
For $1\leq i < j \leq n$ one has $V_j \cap V_i = K \cap K_f$, giving
\begin{equation}\label{sizeorgap}
\left| K_f \setminus \bigcup_{i=1}^{n} K_i \right| = q^{D} - \left( \sum_{j=1}^{n} |V_j| \right) + (n-1)\left| K \cap K_f  \right|.
\end{equation}

Suppose first that $|K \cap K_f| = q^{D-2}$. In this case $K \cap K_f=K$, since $K\cap K_f$ is a subspace of $K$ and is of the same dimension. If each $V_j$ is of size $q^{D-2}$ then $V_j = K$ for each $j$, so $K_f$ is a petal, contradicting the maximality of $n$. Therefore there exists $V_j$ properly containing $K$, so $|V_j| = q^{D-1}$ for some $j$. Now (\ref{sizeorgap}) gives
\[
\left| K_f \setminus \bigcup_{i=1}^{n} K_i \right| \leq q^{D} - q^{D-1},
\] 
since $|V_j| \geq q^{D-2}$ for $1\leq j \leq n$. 

If instead $|K \cap K_f| \leq q^{D-3}$, then (\ref{sizeorgap}) becomes
\[
\left| K_f \setminus \bigcup_{i=1}^{n} K_i \right| \leq q^{D} - nq^{D-2} + (n-1)q^{D-3},
\]
completing the proof.
\end{proofcap}

Examining the proof of Lemma \ref{afterparty}, we see that there are two important subsets of polynomials in $\cF$:
\begin{itemize}
\item a codimension 4 sunflower $\cS = \{ f_1, \ldots, f_n \} \subseteq \cF$ of maximal size, with core $K$,
\item a set $\cS' \subseteq \cF \setminus \cS$ defined by
\begin{equation} \label{Sprime}
\cS' := \left\{ f \in \cF \setminus \cS \colon \ker (f) \cap K = K \right\},
\end{equation}
noting from the proof of Lemma \ref{afterparty} that if $f \in \cS'$ then there exists $i \in \{1,2, \ldots, n\}$ such that $K$ is a proper subspace of $\ker (f) \cap \ker (f_i)$.
\end{itemize}
We will use this notation for the remainder of the section.

\section*{The structure of a sunflower}

As discussed, Lemma \ref{afterparty} suggests that the size of a maximal codimension 4 sunflower plays an important role. In this subsection we establish an upper bound on the size of a sunflower.

\begin{lem}[Maximal size of a sunflower]\label{max_sunfl}
Let $\cS \subseteq \cF$ be a sunflower of codimension $4$. Then
\[
\left| \cS \right| \leq 1 + \frac{q^2+q}{2}.
\]
\end{lem}

We prove the above lemma by exploiting the fact that the sunflowers in $\cF$ come in two types, both easy to bound in size. Before stating that result, we prove that the core has a particular structure, a fact that will be used often in the remainder of this section.

\begin{lem}[The core]\label{core}
Suppose that $\cS$ is a sunflower. Define 
\[
K_{f,f'} = \langle f, Tf, f',Tf' \rangle^\perp.
\]
Then the core $K$ satisfies $K = K_{f,f'}$ for all $f,f' \in \cS$, $f\neq f'$.
\end{lem}
\begin{proofcap}
By the definition of a sunflower, for all distinct $f,f' \in \cS$ we have
\[
K = \ker (f) \cap \ker (f').
\]
Note that $\ker (f) \cap \ker (f')$ corresponds to the kernel of the matrix with rows $f,Tf, f', Tf',$ from which the claim follows.
\end{proofcap}

We are ready to prove that there are only two types of sunflowers. 

\begin{prop}[Structure of a sunflower]\label{structure_sunfl}
For every sunflower $\cS \subseteq \cF$ one of the following holds.
\begin{enumerate}[(i)]
\item (TYPE I) There exists a polynomial $P$ such that for all $f\in \cS$, one can find a polynomial $Q_f$ such that $f = Q_f P$, with $Q_f$ being at most quadratic.
\item (TYPE II) For all $f, f', f'' \in \cS$ we have $f = \lambda f' +  \mu f''$, for some scalars $\lambda, \mu \in \FQ$.
\end{enumerate}
\end{prop}
\begin{proofcap}
Let $f_1 \neq f_2$ be arbitrary elements of $\cS$, and let $P=(f_1, f_2)$. Since each $f \in \cS$ forms a sunflower with $f_1, f_2$, we have $\langle f_1, f_2, Tf_1, Tf_2 \rangle^\perp \subseteq K_f$, and thus $f\in \langle f_1, f_2, Tf_1, Tf_2 \rangle$. Therefore we can find polynomials $L_1, L_2$ of degree at most 1 such that 
\[
f = L_1 f_1 + L_2 f_2 \equiv 0 \pmod{P},
\]
which yields $P\mid f$ for all $f \in \cF$. Thus, defining $g_f = f/P$, note that the set
\[
\cS_P := \{ g_f \colon f \in \cS \}
\]
is a well-defined set of polynomials. Moreover, it constitutes a sunflower, for if $g_{f_1}, g_{f_2} \in \cS_P$ are distinct then
\begin{align*}
(K_{g_{f_1}} \cap K_{g_{f_2}})^\perp &=  \langle f_1/P, f_2/P, Tf_1/P, Tf_2/P  \rangle = \{ h/P : h \in \langle f_1, f_2, Tf_1, Tf_2 \rangle \}  \\
&= \{ h/P: h \in K^\perp \}.
\end{align*}
Finally, note that for all distinct $g,g' \in \cS_P$ we have $(g,g') = 1$.

Now, let $g,g',g'' \in \cS_P$ be pairwise distinct. Since $\cS_P$ is a sunflower, we can find polynomials $L_0, L_1, L_2, L_2'$ of degree at most 1 such that
\begin{align}\label{firstg}
g &= L_1 g' + L_2 g'',  \\
\notag g' &= L_0 g + L_2' g'',
\end{align}
giving
\begin{equation*}
(1- L_0 L_1)g =  (L_1 L_2' + L_2) g''.
\end{equation*}

Suppose first that $L_0 L_1 \neq 1$. Since $g,g'' \in \cS_P$ are coprime, we see that they are at most quadratic, and by a similar argument $\deg(g') \le 2$, so we are in the TYPE I case.

Next, suppose instead that $1 - L_0 L_1 = L_1 L_2' + L_2 = 0$, giving $L_0, L_1 \in \FQ^\times$. Since $\cS_P$ is a sunflower, we can also find polynomials $L_0', L_1'$ of degree at most 1 such that
\[
g'' = L_0' g + L_1' g',
\]
which combined with \eqref{firstg} yields
\begin{equation*}
(1- L_0' L_2)g =  (L_1' L_2 + L_1) g'.
\end{equation*}
If $L_0' L_2 \neq 1$ we are again in the TYPE I case. Otherwise $L_0', L_2 \in \FQ^\times$, and by symmetry we are in the TYPE II case.
\end{proofcap}

We are ready to bound the size of a sunflower.

\begin{proofof}{Lemma \ref{max_sunfl}}
Recall that we assumed without loss that all $f \in \cF$ are monic. If $\cS$ is of TYPE II, then counting the polynomials projectively yields \[|\cS| \leq \left| \PP_{\FQ}^1 \right| = \frac{q^2-1}{q-1} = q+1.\]

Now suppose $\cS$ is of TYPE I, with $P$ defined accordingly. We may assume that $P$ is monic. Let $K$ be the core of $\cS$. Then $K$ contains the codimension 4 subspace $\{ PC: \deg \: C \le 3 \}^\perp$, and so
\[
K = \{ PC: \deg \: C \le 3 \}^\perp.
\]
For each pair $f \neq f'\in \cS$, we now have
\begin{equation}\label{dimensionzero}
\langle f,Tf, f', Tf' \rangle = \{ PC: \deg \: C \le 3 \}.
\end{equation}

Note that $\{ f/P: f \in \cS \}$ is a set of monic polynomials of degree at most 2. Fix an ordering $\lambda_1, \ldots, \lambda_q$ of $\FQ$ and consider the following colouring of the set of monic polynomials of degree at most 2.
\begin{itemize}
\item The polynomial 1 has its own colour $c$;
\item each irreducible quadratic $Q$ has its own colour $c_Q$;
\item for $1\leq i \leq q$, multiples of $T-\lambda_i$, but not of $T-\lambda_j$ for $j<i$, have colour $c_i$.
\end{itemize}
This induces a colouring of $\cS$ in the obvious way. Observe from (\ref{dimensionzero}) that distinct $f$ and $f'$ satisfy $\gcd(f,f') = P$. Therefore the sunflower cannot contain two polynomials of the same colour, so the total number satisfies
\[
|\cS| \leq 1 + \frac{q^2-q}{2} + q,
\]
since the number of monic irreducible quadratics is $\frac12 (q^2-q)$.
\end{proofof}


\section*{The case when there exists a large sunflower}

Here we consider the case in which there is a large codimension 4 sunflower in the set of speeds, that is, a sunflower of size $n\geq q+1$. In this case, we see from Lemma \ref{afterparty} that all but $n$ polynomials cover at most $q^D - q^{D-1}$ new points. This gives a non-trivial bound on the number of polynomials needed. 

\begin{prop}[Large sunflower implies halfway]\label{largeimpliesdone}
Let $\cF$ be of loneliness at most $q^{-3}$. If there exists a sunflower of codimension $4$ and size at least $q+1$, then
\[
|\cF| \geq q^2 + \frac{q+1}{2}.
\]
\end{prop}

\begin{remark}\label{remark_large}
In the large-sunflower case, one can work harder to establish the sharp inequality $|\cF| \geq q^2 + q+1$. To ease the exposition, we presently only prove a weaker result that corresponds to the bound we obtain in the case when all sunflowers are small. We provide, in the appendix, a proof of the inequality $|\cF| \geq q^2 + q+1$, assuming the existence of a sunflower of codimension 4 and size at least $q+2$, and that $q \ge 8$.
\end{remark}

\begin{proofcap}
By Lemma \ref{lrc_covering} we see that $\bigcup_{f\in\cF} \ker (f) = \FQ^{D+2}$. Let $\cS$ be a codimension 4 sunflower of maximal size, and let $n= |\cS|$. Since $n \geq q+1$, one has 
\[
\max \{ q^D - q^{D-1}, \hspace{3mm}q^{D} - n q^{D-2} + (n-1) q^{D-3} \} = q^D - q^{D-1}.
\]
Defining $R = |\cF|$ and using Lemma \ref{afterparty}, we get
\begin{align*}
q^{D+2} &= \left| \FQ^{D+2} \right| = \left| \bigcup_{f\in\cF} \ker (f) \right| \leq \left| \bigcup_{f'\in\cS} \ker (f')  \right| + \sum_{f\in \cF \setminus \cS} \left| \ker (f) \setminus \bigcup_{f'\in\cS} \ker (f') \right| \\
&\leq \left( q^D + (n-1)(q^D - q^{D-2}) \right) + (R - n)(q^D - q^{D-1}). 
\end{align*}
Expanding the right-hand side gives
\[
R \geq q^2 + q + \frac{q+1-n}{q},
\]
which together with Lemma \ref{max_sunfl} yields
\[
R \geq q^2 + \frac{q +1}{2}.
\]
\end{proofcap}


\section*{The case when all sunflowers are small}

Here we discuss the remaining case, where all codimension 4 sunflowers have at most $q$ petals. This turns out to be much more demanding, as the bound \eqref{after_contribution} becomes less powerful. In what follows we are able to improve on (\ref{after_contribution}) in this case, by showing that there are many distinct pairwise intersections. In doing so we have to increase the influence of three-fold intersections, and thus are not able to establish Conjecture \ref{thm:lrcffconj} in full. However, the aforementioned improvement of (\ref{after_contribution}) yields a significant gain over Theorem \ref{maximalloneliness} in the case $k=2$.

\begin{prop}[Small sunflowers imply almost halfway]\label{small_done}
Let $\cF$ be such that all sunflowers in $\cF$ of codimension $4$ are of size at most $q$. If $\cF$ is of loneliness at most $q^{-3}$, then there exists a universal constant $C \in \bR$ such that $|\cF| \geq q^2 + 0.4877q - C$.
\end{prop}

The rest of the paper will be devoted to the proof of this proposition. By Lemma \ref{lrc_covering}, we have that loneliness at most $q^{-3}$ implies $\bigcup_{f\in\cF} \ker (f) = \FQ^{D+2}$. Define $R = |\cF|$, and let $n \le q$ be the size of a maximal codimension 4 sunflower $\cS$ in $\cF$. We order the runners $f_1, f_2, \ldots, f_R$ in such a way that the two-fold intersections all have codimension 4, until the \textit{change point}, after which all the contributions are at most $q^{D} - q^{D-1}$.

As in the proof of Proposition \ref{largeimpliesdone}, we bound the contribution from the first $n$ runners (those that are in $\cS$) by a total of $q^{D} + (n-1)(q^{D}-q^{D-2})$, and the next $(n-1)^2$ runners by $q^{D} - n q^{D-2}+(n-1)q^{D-3}$ each, using Lemma \ref{afterparty}. We call this the \textit{initial phase}.

Next we consider runner $m > n(n-1) +1$. Using the same notation as in the proof of Lemma \ref{afterparty}, namely that $K_i = \ker (f_i)$ and $V_i = K_i \cap K_m$, the covering contribution of $f_m$ is 
\[
q^D - |V_1 \cup \cdots \cup V_{m-1}|.
\]
As mentioned in the paragraph at the beginning of this case, the main idea is that there exist sufficiently many distinct two-fold intersections. 

\begin{claim}\label{claim_one}
Let $t \ge 2$ be an integer, and assume that at $m$ we have not yet reached the change point. Then there are at least $t$ distinct sets among $V_1, \ldots, V_{m-1}$, provided that
\begin{equation}\label{t_distinct}
m-1 > (t-1)(n-1).
\end{equation}
Moreover, the sets $V_1, \ldots, V_n$ are distinct, and the core $K$ of $\cS$ satisfies 
\[
| K \cap K_m | \leq q^{D-3}.
\]
\end{claim}

\begin{proofcap}
Since we have not passed the change point, all of the $V_i$ have codimension $4$. If $V_i = V_j$, then $K_i,K_j,K_m$ form a sunflower of codimension 4, as $K_i \cap K_j \subseteq V_i = V_j$. Since $n$ is the maximum size of a sunflower, we see that for each $V_i$ there can be at most $n-1$ indices $j \in \{1,2,\ldots,m-1\}$ such that $V_i = V_j$. In view of the inequality \eqref{t_distinct}, the pigeonhole principle ensures that there are at least $t$ distinct sets among $V_1, \ldots, V_{m-1}$.

For the second part, assume for a contradiction that for some $1\leq i < j \leq n$ we have $V_i = V_j$. Then the core $K$ of $\cS$ satisfies 
\[
K\cap K_m = V_i \cap V_j = V_i = V_j,
\]
implying that $K \cap K_m$ has codimension $4$. Now $K \cap K_m = K$, so by the definition of $\cS'$ following Lemma \ref{afterparty}, we have that $f_m \in \cS'$. However, as discussed following that definition, this implies that there exists $i \in \{1,2,\ldots,n \}$ such that $K$ is a proper subspace of $V_i$, contradicting the assumption that we have not passed the change point. Therefore the sets $V_1, \ldots, V_n$ are pairwise distinct, and so $ |K\cap K_m | \leq q^{D-3}$.
\end{proofcap}

\begin{claim}\label{claim_two}
The next $(n-1)(q-n)$ runners after the initial phase contribute a total of at most 
\[
(n-1) \sum_{t=n+1}^{q} \left( q^{D} - tq^{D-2} + \left[ {t \choose 2}  - {n-1 \choose 2} \right] q^{D-3} \right).
\]
\end{claim}
\begin{proofcap}
Consider runner $m > n(n-1) + 1$, and suppose for the time being that this process precedes the change point. If $m$ is among the first $n-1$ runners after the initial phase we may assume, by Claim \ref{claim_one}, that $V_1, \ldots, V_{t}$ are distinct with $t=n+1$. In this case the contribution of runner $m$ is at most
\begin{equation}\label{contribution_t_after}
q^{D} - \sum_{j=1}^{t} |V_j| + (n-1) |K\cap K_m| + \sum_{j=n+1}^{t} \left|V_j \cap (V_1 \cup \cdots \cup V_{j-1})  \right|,
\end{equation}
since $V_i \cap V_j = K \cap K_m$ for all $1\leq i < j \leq n$. Note that for distinct $V_i$ and $V_j$ the codimension of $V_i \cap V_j$ is at least 5, so the quantity in \eqref{contribution_t_after} is bounded from above by
\begin{align*}\label{contribution_t}
& q^{D} - tq^{D-2} + (n-1) q^{D-3} + q^{D-3} \sum_{j=n+1}^{t} (j-1) 
\\ &\qquad =  q^{D} - tq^{D-2} + \left[ {t \choose 2} - {n-1 \choose 2} \right] q^{D-3}.
\end{align*}

The contribution from the next $n-1$ runners is the same, but with $t=n+2$, and so on. To conclude, observe that the result remains valid even if we cross the change point, as the above bounds on the individual contributions exceed $q^{D} - q^{D-1}$.
\end{proofcap}

Now suppose we are in the \textit{final stage}, meaning that
\[
m > q(n-1) +1.
\]
Then we may apply Claim \ref{claim_one} with $t=q$. The displayed expression in Claim 2 has a minimum at $t=q+\frac12$, so there is nothing to be gained from choosing a larger value of $t$. Therefore the contribution of each runner in the final stage is bounded by 
\begin{equation}\label{final_stage}
q^{D} - \frac{q^{D-1} + q^{D-2}}{2} - {n-1 \choose 2}q^{D-3},
\end{equation}
which holds irrespective of whether or not we have crossed the change point. 

Combining everything, we find that
\begin{equation}\label{final_calculation}
q^{D+2} = \left| \FQ^{D+2} \right| = \left| \bigcup_{i=1}^{R} K_i \right| \leq S_1 + S_2 + S_3,
\end{equation}
where 
\[
S_1 = q^{D} + (n-1) (q^D - q^{D-2}) + (n-1)^2 (q^D - nq^{D-2} + (n-1)q^{D-3})
\]
comes from the initial phase,
\[
S_2 = (n-1) \sum_{t=n+1}^{q} \left( q^{D} - tq^{D-2} + \left[ {t \choose 2} - {n-1 \choose 2} \right]q^{D-3} \right)
\]
comes from Claim \ref{claim_two}, and
\[
S_3 = (R-1-q(n-1)) \left( q^{D} - \frac{q^{D-1} + q^{D-2}}{2} - {n-1 \choose 2}q^{D-3}\right)
\]
comes from the final stage, using (\ref{final_stage}) for these final polynomials.

We may assume without loss that $R = O(q^2)$. Simplifying the right-hand side of (\ref{final_calculation}), dividing by $q^{D-3}$, and then grouping together all the summands that are $O(q^3)$, gives
\begin{align*}
& R \left( q^3 - \frac{q^{2} + n^2}{2}\right) \\
 &\qquad \ge q^{5} + n^3q  - \frac12 nq^3 - \frac12 n^4 + n \sum_{t=n+1}^{q} \left( tq - {t \choose 2}  \right) - C_1 q^3,
\end{align*}
since $n \leq q$. Note that
\begin{align*}
\left( q+\frac12 \right) \sum_{t=n+1}^{q} t - \frac12 \sum_{t=n+1}^{q} t^2 = q\frac{q^2-n^2}{2} + \frac{n^3 - q^3}{6} +  O(q^2),
\end{align*}
which yields
\begin{align*}
&R\left(q^3 - \frac{q^2+n^2}{2} \right) \\
&\qquad \geq q^2 \left(q^3 - \frac{q^2+n^2}{2} \right)  - \frac{2n^4 - 3n^3q -3n^2q^2+ nq^3 - 3q^4}{6} - C_2 q^3.
\end{align*}

Determining $\lambda = n/q \in [0,1]$ numerically to maximise 
\[
\frac{2n^4 - 3n^3q -3n^2q^2+ nq^3}{q^4} = 2 \lambda^4 - 3 \lambda^3 - 3 \lambda^2 + \lambda,
\]
we see that
\[
R\left(q^3 - \frac{q^2+n^2}{2} \right) \geq q^2 \left(q^3 - \frac{q^2+n^2}{2} \right) + 0.4877q^4 - C_2 q^3,
\]
giving
\[
R \geq q^2 + 0.4877q -C.
\]
This completes the proof of Proposition \ref{small_done}.

Finally, Propositions \ref{largeimpliesdone} and \ref{small_done} imply Theorem \ref{close_half}.

\appendix
\section{Large sunflowers}

Here we establish the following strong form of Proposition \ref{largeimpliesdone}, as promised in Remark \ref{remark_large}.

\begin{prop}\label{large_main_prop_done}
Let $\cF$ be of loneliness at most $q^{-3}$. If $q \ge 8$ and there exists a sunflower of codimension $4$ and size at least $q+2$, then
\[
|\cF| \geq q^2 + q + 1.
\]
\end{prop}

We start by recalling the proof of Lemma~\ref{afterparty} and, in particular, the definitions of $\cS$ and $\cS'$, which we will use for the rest of this appendix. Let $\cS$ be a sunflower of codimension 4 that is of maximal size in $\cF$, denote its core by $K$, and define $\cS'$ by \eqref{Sprime}. Put \[
n:= |\cS|, \hspace{5mm}k:= |\cS'|,\]
observing that $n\geq q+ 2$. We order the polynomials in $\cF$ so that \[\cS =\{f_1,\ldots,f_n\}, \qquad \cS' = \{ f_{n+1}, \ldots, f_{n+k} \},\]
and define $\cS'' := \cF \setminus (\cS \cup \cS')$. Lemma~\ref{afterparty} and its proof imply that for each $f\in \cS'$ we have
\begin{equation*}\label{cS1_contribution}
\left| \ker (f) \setminus \bigcup_{i=1}^{n} \ker (f_i) \right| \leq q^{D} - q^{D-1},
\end{equation*}
and for each $f\in \cS''$ we have
\begin{equation}\label{cS2_contribution}
\left| \ker (f) \setminus \bigcup_{i=1}^{n} \ker (f_i) \right| \leq q^{D} - nq^{D-2}+(n-1)q^{D-3} < q^{D} - q^{D-1}.
\end{equation}

Our strategy is as follows. As before, we consider the covering problem arising from Lemma~\ref{lrc_covering}. We first show that in order for $\displaystyle \bigcup_{f\in\cF} \ker (f) = \FQ^{D+2}$ to hold when there exists a large sunflower in the set of speeds, one needs both $n$ and $k$ to be close to maximal. Moreover, we prove that $\cS$ and $\cS'$ are ``closely connected": loosely speaking, the sunflower $\cS$ should be thought of as comprising almost all irreducible polynomials of degree at most $2$, and $\cS'$ comprising almost all reducible quadratic polynomials. Then $\cS''$ has to include any remaining quadratics; thus, $\cF$ has to look like the example in (\ref{counterexample}), up to multiplication by a polynomial.

\begin{lemma}[$\cS \cup \cS'$ is large]\label{large_structured_part}
Let $\cF$ be of loneliness at most $q^{-3}$ and let $\cS$ be a sunflower of codimension 4 that is of maximal size. If $|\cF| \leq q^2+q$ and $|\cS|\geq q+2$, then for $\cS'$ defined as above one has $|\cS|+|\cS'| \geq q^2.$ 
\end{lemma}
\begin{proofcap}
Recall the definitions of $n$ and $k$, and let $r_n, r_k \in \bR$ be such that \[n= |\cS| = \frac12 q^2 + r_n,  \hspace{5mm}k = |\cS'| =\frac12 q^2 + r_k.\] 
Then the assumption that $\FF$ is of loneliness at most $q^{-3}$, combined with Lemma \ref{lrc_covering} and the above discussion, yields
\begin{align}\label{rnrk_bound1}
q^{D+2} &\leq \underbrace{\left( q^D + \left(\frac{q^2}{2}+r_n -1\right)(q^D - q^{D-2})\right)}_{=S_1} + \underbrace{\left(\frac{q^2}{2} + r_k \right)(q^D - q^{D-1})}_{=S_2}  \nonumber \\
&\hspace{20mm}+ \underbrace{(q-r_n - r_k) \left(\frac12 q^D -  r_n q^{D-2} + \frac12 q^{D-1} + (r_n-1)q^{D-3}\right)}_{=S_3}.
\end{align}
Here $S_1, S_2$ and $S_3$ are upper bounds for the covering contributions of $\cS$, $\cS'$ and $\cS''$, respectively. We compute that \begin{align}\label{rnrk_bound}
S_1 + S_2 + S_3 &= q^{D+2} + (r_n+r_k)q^{D-3} \left( \frac12q^3 -\frac32 q^2 + r_n q - (r_n-1)\right) \nonumber \\
&= q^{D+2} + (r_n+r_k)q^{D-3} (q-1)(n-(q+1)). 
\end{align}
Since $n > q+1$, we deduce from (\ref{rnrk_bound1}) and (\ref{rnrk_bound}) that $r_n+r_k \geq 0$, completing the proof.
\end{proofcap}

In the course of the proof of Lemma \ref{max_sunfl}, we showed that any TYPE II sunflower (in the sense of Proposition \ref{structure_sunfl}) has size at most $q+1$. As $n>q+1$, the sunflower $\cS$ must therefore be of TYPE I. In other words, there exists a polynomial $P$ and at most quadratic polynomials $Q_1,\ldots, Q_n$ such that $f_i = Q_iP$, for all $i \in \{1,2,\ldots, n\}$. 
We now show that $\cS'$ has the same structure, with the additional information that most of the polynomials involved are reducible quadratics.

\begin{prop}[$\cS'$ is an extension of $\cS$]\label{s1_structure}
For each $f \in \cS'$, there exists a polynomial $Q_f$, at most quadratic, such that $f = PQ_f$. Moreover, if $Q_f$ is of degree $2$, then $Q_f$ is reducible.
\end{prop}
\begin{proofcap}
First we prove that $f = PQ_f$ for some polynomial $Q_f$ of degree at most 2. The way that $\cS'$ is constructed implies that $K \cap \ker (f) = K$. Hence we can find $L_{1},L_1',L_2,L_2'$, at most linear, such that
\begin{align*}
f &= L_1f_1 + L_2 f_2 = P(L_1 Q_{1} + L_2 Q_{2}),\\
Tf &= L_1'f_1 + L_2'f_2 = P(L_1' Q_{1} + L_2' Q_{2}).
\end{align*}
The first equation gives $P \mid f$, and the second implies that $\deg (f) \leq \deg (P) +2$, which combine to give the existence of an at most quadratic $Q_f$ such that $f = PQ_f$. 

Now suppose that $Q_f$ is of degree 2. Since $f \in \cS'$, we can find $ i\in \{1,2,\ldots, n \}$ such that $|\ker (f)  \cap \ker (f_i) | = q^{D-1}$. Therefore, there exist polynomials $L$ and $L_{i}$, at most linear, such that 
\[
L PQ_f = L f  = L_{i} f_i = L_{i} P Q_i.
\]
Observe that $Q_f$ being an irreducible quadratic would imply $f = f_i$. Therefore $Q_f$ must instead be reducible. 
\end{proofcap}

As an immediate consequence we obtain $|\cS'| \leq \frac12 (q^2+3q+2)$. The aforementioned upper bounds on $|\cS|$ and $|\cS'|$, together with Lemma \ref{large_structured_part}, imply that both $|\cS|$ and $|\cS'|$ are close to maximal. Moreover, Propositions \ref{structure_sunfl}, \ref{s1_structure}, and Lemma \ref{large_structured_part} imply that at least $q^2$ polynomials in $\cF$ are of the form $f=PQ_f$, with $Q_f$ at most quadratic. This brings us very close to showing that $\cF$ indeed looks like the example in (\ref{counterexample}). 

The fact that $\cS$ and $\cS'$ make up most of the set $\cF$ allows us to improve on the bound in Lemma \ref{afterparty}. Informally, the following lemma asserts that it is sufficient to be able to add a small set as a ``bridge" between $\cS$ and $\cS'$, in order that for the remaining polynomials the contribution bound in Lemma \ref{afterparty} may be improved to $q^{D}-2q^{D-1}+q^{D-2}$.

\begin{lemma}[Sufficient connections between $\cS$ and $\cS'$]\label{s_s1_connections}
Let $\FF$ be of loneliness at most $q^{-3}$ such that a sunflower $\cS \subseteq \cF$ of maximal size satisfies $|\cS| \geq 2q+1$. Suppose that there exists $\cC \subseteq \FF \setminus \cS$ such that 
\begin{equation}\label{S_C_size}
2|\cS| +  |\cC| < (q+1)^2,
\end{equation} 
and such that for all $f\in \cS' \setminus \cC$ there exist distinct $f_i, f_j \in \cS \cup \cC$ satisfying 
\begin{equation}\label{s_c_ker}
|\ker (f) \cap \ker (f_i) | = |\ker (f) \cap \ker (f_j)| = q^{D-1}, \hspace{3mm} |\ker (f_i) \cap \ker (f_j) | = q^{D-2}.
\end{equation}
Then
\[
|\FF| \geq q^2+q+1.
\]
\end{lemma}

\begin{remark} \label{nothing} The assumption on $|\cS|$ in Lemma \ref{s_s1_connections} is stronger than in previous instances of the large-sunflower case. As $q \ge 8$, this will cost us nothing, since Lemma \ref{s_s1_connections} will only be applied when $|\cS| \approx q^2/2$.
\end{remark}

\begin{proofcap}
Let $R:=|\FF|$ and $r:=|\cC|$. Lemma \ref{afterparty} implies that the covering contribution of each $f\in \cC$ is at most $q^{D}-q^{D-1}$, since $\cC \subseteq \cS' \cup \cS''$ and $n = |\cS| \geq 2q+1$. 

Let $f\in \cS' \setminus \cC$, and let $f_i, f_j \in \cS \cup \cC$ be such that (\ref{s_c_ker}) holds. Then the contribution of $f$ is bounded from above by
\begin{align*}
\left|\ker (f) \setminus \left( \ker (f_i) \cup \ker (f_j) \right) \right| \leq q^D - 2q^{D-1} + q^{D-2}.
\end{align*}
Moreover, we know from \eqref{cS2_contribution} that the contribution of any $f \in \cS''$ is at most
\[
q^{D}-nq^{D-2} + (n-1)q^{D-3} \leq q^{D-2}(q-1)^2,
\]
since $n\geq 2q+1$. Thus, for all $f\in \cS' \cup \cS''  \setminus \cC$ the contribution of $f$ is at most $q^{D-2}(q-1)^2$. Summing all the contributions and using the assumption that $\del(\FF) \le q^{-3}$ via Lemma \ref{lrc_covering}, one has
\begin{align*}
q^{D+2} &\leq q^D + (n-1)(q^D - q^{D-2}) + r(q^D - q^{D-1}) + (R-n - r)q^{D-2}(q-1)^2.
\end{align*} 
This in turn implies that
\[
R \geq q^2 + q + \frac{(q+1)^2 - (2n+r)}{q-1}> q^2+q,
\]
using the assumption $2n+r < (q+1)^2$.
\end{proofcap}

\bigskip
We are ready to establish the main result of this section. The key step is to construct the set $\cC$ appearing in Lemma \ref{s_s1_connections}.

\begin{proofof}{Proposition \ref{large_main_prop_done}}
Let $R=|\FF|$, $n = |\cS|$ and $k= |\cS'|$, as before. For the sake of obtaining a contradiction, suppose $R \leq q^2+q$. By increasing the size of $\cF$ if needed, we may in fact assume that 
\begin{equation} \label{wma}
R = q^2 + q.
\end{equation}
Recall from Lemma \ref{large_structured_part} that
\begin{equation}\label{n_klarge}
n+k \geq q^2.
\end{equation}
As $q \ge 8$ we have
\[
k \le \frac{q^2+3q+2}2 \le q^2 - 2q - 1,
\]
so \eqref{n_klarge} ensures that $n \ge 2q + 1$. This formalises Remark \ref{nothing}, and confirms one of the hypotheses of Lemma \ref{s_s1_connections}.\\

\noindent \textbf{Claim.} For each $\lambda \in \FQ$ we can find an at most linear polynomial $L$ such that \[P(T-\lambda)L \in \cF.\] 

\begin{proofcap}
Suppose to the contrary that there exists $\lambda \in \FQ$ such that 
\[
P(T-\lambda)L \notin \FF
\]
for all $L$ that are at most linear. Then (\ref{n_klarge}), together with Propositions \ref{structure_sunfl} and \ref{s1_structure}, implies that
\[
\cS \cup \cS' = \left\{ PQ \colon (T-\lambda) \nmid Q\right\},
\]
since $\left|  \left\{ PQ \colon (T-\lambda) \nmid Q\right\} \right| = q^2$. Hence the largest sunflower has $\frac12 (q^2+q)$ petals, and is formed by multiplying a polynomial $P$ by all (monic) irreducible quadratic polynomials, and all $L^2$, for $L\neq T-\lambda$ at most linear. We wish to prove that in this case $\cS$ satisfies the assumptions of Lemma \ref{s_s1_connections}, with $\cC = \emptyset$. Indeed, let $f = PL_i L_j \in \cS'$ be arbitrary. Since $L_i,L_j \neq T-\lambda$, it is easy to see that (\ref{s_c_ker}) holds with $f_i := PL_{i}^2$ and $f_j := PL_{j}^2 \in \cS$. Now Lemma \ref{s_s1_connections} delivers a contradiction, as we assumed $R \leq q^2+q$.
\end{proofcap}

\bigskip

 Let
\[
\cS_{\max}:= \left\{ \cS \subset \FF \colon \cS\text{ is a maximal sunflower} \right\}.
\]
To each $\cS \in \cS_{\max}$ we associate the set
\[
\lambda_{\cS} := \left\{ \lambda \in \FQ \colon \exists Q_\lambda \text{ such that } (T-\lambda)|Q_\lambda \text{ and }PQ_\lambda \in \cS \right\}.
\]
The colouring argument in the proof of Lemma \ref{max_sunfl} ensures that for $\lambda \in \lambda_{\cS}$ the polynomial $Q_\lambda$ is unique. Define $\lambda_{\cS}^c := \FQ \setminus \lambda_{\cS}$, and let $\cS \in \cS_{\max}$ be such that $|\lambda_{\cS}^c|$ is maximal, and such that if $PL^2 \in \cF$ then $PL \notin \cS$. 


Observe that if $\lambda, \lambda' \in \lambda_{\cS}^c$ then $P(T-\lambda)(T-\lambda') \notin \FF$, for otherwise we could add $P(T-\lambda)(T-\lambda')$ to $\cS$, contradicting its maximality. In particular $P(T-\lambda)^2 \notin \FF$ for $\lambda \in \lambda_{\cS}^c$. With $\ell := |\lambda_{\cS}^c|$ we have
\[
n \leq 1 + \frac{q^2+q}{2} - \ell;
\]
this bound follows from variant of the colouring argument in the proof of Lemma \ref{max_sunfl}.

Roughly speaking, we want $\cS$ to attain the fewest distinct elements of $\FQ$ as roots, so that we may construct $\cC$ efficiently. Choosing $\cS$ as above ensures that if $L$ is at most linear then $PL^2 \notin \cF \setminus \cS$. Indeed, first consider when $L=1$, supposing for a contradiction that $P \in \cF \setminus \cS$. If $\cS$ contains no linear multiples of $P$, then we can contradict its maximality by including $P$. If, on the other hand, $\cS$ contains a linear multiple of $P$, then we can put $P$ in $\cS$ instead of it, to contradict the maximality of $|\lam_S^c|$. We conclude that $P \notin \cF \setminus \cS$. 

Otherwise $L= T - \lam$. Suppose for a contradiction that $PL^2 \in \cF \setminus \cS$. Then we must have $\lam \in \lam_S$, so $PLL' \in \cS$ for some $L'$ at most linear. Moreover, we must have $L' = 1$, for otherwise we could replace $PLL'$ by $PL^2$ in $\cS$ to increase $|\lam_\cS^c|$. Now $PL^2 \in \cF$ and $PL \in \cS$, which is impossible, by our choice of $\cS$.\\

We split our argument into two cases. Case A applies when for each $\lambda \in \lambda_{\cS}^c$ we can find two distinct polynomials $L_{1}^{(\lambda)}, L_{2}^{(\lambda)}$, of degree at most 1, such that 
\begin{equation}\label{f_lambda12}
f_{\lambda,1} = P(T-\lambda)L_{1}^{(\lambda)}, \hspace{5mm}f_{\lambda,2} = P(T-\lambda)L_{2}^{(\lambda)} \in \FF.
\end{equation}
Case B applies when there exists $\lambda \in \FQ$ for which there is only one \[f_\lambda = P(T-\lambda)L_\lambda \in \FF.\]

\bigskip
First consider Case A. For $\lambda \in \lambda_{\cS}^c$ choose distinct $f_{\lambda,1}, f_{\lambda,2} \in \FF$ as in (\ref{f_lambda12}) and add them to $\cC$, with $L_{\lambda,1} = 1$ if $P(T-\lambda)\in\FF$. Note that $\cS$ and $\cC$ are disjoint, and
\[
2 |\cS| + |\cC| \leq 2\left(1+\frac{q^2+q}{2} - \ell \right) + 2 \ell < (q+1)^2.
\]

We proceed to confirm the hypotheses (\ref{s_c_ker}) of Lemma \ref{s_s1_connections}. Let 
\[
f = P(T-\lambda_i)(T-\lambda_j) \in \cS' \setminus \cC,
\]
and observe that $\lambda_i \neq \lambda_j$, since we showed that if $PL^2 \in \cF$ then $PL^2 \in \cS$. If $\lambda_i,\lambda_j \in \lambda_{\cS}$ then we can find $f_i = P(T-\lambda_i)L_i$ and $f_j = P(T-\lambda_j)L_j$ in $\cS$, whereupon (\ref{s_c_ker}) holds.
If $\lambda_i \in \lambda_{\cS}$ and $\lambda_j  \in \lambda_{\cS}^c$, then we can find $f_i = P(T-\lambda_i)L_i \in \cS$ and \[f_{\lambda_j,1} = P(T-\lambda_j)L_{1}^{(\lambda_j)}, \hspace{3mm} f_{\lambda_j,2} = P(T-\lambda_j)L_{2}^{(\lambda_j)} \in \cC.\]
Choose $m \in \{1,2\}$ such that $L_{m}^{(\lambda_j)} \neq L_i$ and define $f_j := f_{\lambda_j,m}$ to get (\ref{s_c_ker}) again, noting that $f \in \cS' \setminus \cC$ implies $f_j \neq f$. Since we proved that there is no polynomial in $\FF$ for which $\lambda_i,\lambda_j \in \lambda_{\cS}^c$, this analysis covers all possibilities. 

To conclude Case A it remains to consider $f = P(T-\lambda_i) \in \cS' \setminus \cC$, and show that it satisfies (\ref{s_c_ker}) for some $f_i, f_j \in \cS \cup \cC$. As $L_{1}^{(\lambda_i)}=1$ had priority over genuinely linear polynomials in the construction of $\cC$, we see that $P(T-\lambda_i) \in \cS' \setminus \cC$ only if $\lambda_i \in \lambda_{\cS}$. Hence there exists $f_i = P(T-\lambda_i)L_i \in \cS$. If there were no linear $L$ such that $PL \in \cS$, then we could have chosen $P(T-\lambda_i)$ to place in $\cS$ instead of $f_i$, and this contradicts the minimality of $|\lambda_{\cS}|$. Therefore there exists $L_j \neq T-\lambda_i$, at most linear, such that $f_j = PL_j \in \cS$, giving (\ref{s_c_ker}).

Putting everything together, we see that $\cS$ and $\cC$ satisfy the assumptions of Lemma \ref{s_s1_connections}, which implies $R\geq q^2+q+1$ as desired. \\

Finally, suppose that we are in Case B, so that there exists $\lambda \in \FQ$ such that there is a unique $f_\lambda = P(T-\lambda)L_\lambda$ in $\FF$. Then, by \eqref{wma}, there is at most one polynomial from the set
\[
X:=  \{ PQ \colon (T-\lambda)\nmid Q, \: \deg \: Q \le 2 \}
\]
that is not in $\FF$, as $|X| = q^2$. We initially choose 
\[
\cC = \begin{cases}
 \{ f_\lam \} \setminus \cS, &\text{if } P \in \cF \\
\{ f_\lam \} \cup \{ P(T-\kap) \in \cF : \kap \in \bF_q \} \setminus \cS, &\text{if } P \notin \cF.
\end{cases}
\]
If $P(T-\lam')^2 \notin \cF$ for some $\lam' \ne \lam$, then for some $\lam'' \notin \{ \lam, \lam' \}$ we replace $P(T-\lam'')^2$ by $P(T-\lam')(T-\lam'')$ in $\cS$, and also append $P(T-\lam'')^2$ to $\cC$. 
When $P \notin \cF$, we note that
\[
2|\cS| + |\cC| \le (q^2 + q) + q < (q+1)^2.
\]
The inequality \eqref{S_C_size} also holds when $P \in \cF$.

We showed in Case A that if $L$ is at most linear then $ PL^2 \notin \cF \setminus \cS$; the argument still works in Case B---wherein $\cS$ has possibly been modified---unless $L = T-\lam''$, and in the latter case $PL^2 \in \cC$. Thus, if $f \in \cS' \setminus \cC $ then either
\begin{enumerate} [(i)]
\item $f = P (T- \lam_i)(T- \lam_j)$ with $\lam_i, \lam_j, \lam \in \bF_q$ pairwise distinct, or
\item $f = P(T-\lam_j)$ for some $\lam_j \ne \lam$.
\end{enumerate}
In scenario (i) we obtain \eqref{s_c_ker} by choosing $f_i = P(T-\lam_i)^2$ and $f_j = P(T- \lam_j)^2$, or instead $f_j = P(T-\lam')(T-\lam'')$ if say $\lam_j \in \{ \lam', \lam'' \}$. In scenario (ii) we must have $P \in \cF$, by our choice of $\cC$, and so $P \in \cS$. We may therefore choose $f_ i = P$ and $f_j = P(T-\lam_j)^2$, as long as $\lam_j \notin \{ \lam', \lam'' \}$. If $\lam_j \in \{ \lam', \lam'' \}$ then we may choose $f_j = P(T-\lam')(T-\lam'')$ instead.

Therefore $\cS$ and $\cC$ satisfy the hypotheses of Lemma \ref{s_s1_connections}, which finishes the proof of Proposition \ref{large_main_prop_done}.
\end{proofof}

\bigskip

\end{document}